\def\draft{n}
\documentclass{amsart}
\usepackage[headings]{fullpage}
\usepackage{amssymb,epic,eepic,epsfig,amsbsy,amsmath}

\theoremstyle{plain}

\newtheorem{theorem}{Theorem}
\newtheorem{proposition}{Proposition}[section]
\newtheorem{lemma}[proposition]{Lemma}
\newtheorem{corollary}[proposition]{Corollary}

\theoremstyle{definition}

\theoremstyle{remark}

\newtheorem{remark}[proposition]{Remark}

\def\printname#1{
    \if\draft y
        \smash{\makebox[0pt]{\hspace{-0.5in}
            \raisebox{8pt}{\tt\tiny #1}}}
    \fi
}

\newlength{\standardunitlength}
\catcode`\@=11
\long\def\@makecaption#1#2{%
     \vskip 10pt

\setbox\@tempboxa\hbox{
       \small\sf{\bfcaptionfont #1. }\ignorespaces #2}%
     \ifdim \wd\@tempboxa >\captionwidth {%
         \rightskip=\@captionmargin\leftskip=\@captionmargin
         \unhbox\@tempboxa\par}%
       \else
         \hbox to\hsize{\hfil\box\@tempboxa\hfil}%
     \fi}
\font\bfcaptionfont=cmssbx10 scaled \magstephalf
\newdimen\@captionmargin\@captionmargin=2\parindent
\newdimen\captionwidth\captionwidth=\hsize
\catcode`\@=12

\def\lbl#1{\label{#1}\printname{#1}}

\def\BZ{\mathbb Z}

\def\BQ{\mathbb Q}

\def\BC{\mathbb C}

\def\D{\Delta}

\def\a{\alpha}

\def\la{\langle}
\def\ra{\rangle}

\def\e{\epsilon}

\def\mat#1#2#3#4{\left(
\begin{matrix}
 #1 & #2  \\
 #3 & #4
\end{matrix}
\right)}

\def\longto{\longrightarrow}

\def\pt{\partial}

\def\SL{\mathrm{SL}}
\def\vol{\mathrm{vol}}
\renewcommand{\Re}{\mathrm{Re}}
\def\coeff{\mathrm{coeff}}
\def\Taylor{\mathrm{Taylor}}

\begin{document}


\title[An analytic version of the MMR Conjecture]{
An analytic version of the Melvin-Morton-Rozansky Conjecture}

\author{Stavros Garoufalidis}
\address{School of Mathematics \\
         Georgia Institute of Technology \\
         Atlanta, GA 30332-0160, USA \\
         {\tt http://www.math.gatech} \newline {\tt .edu/$\sim$stavros } }
\email{stavros@math.gatech.edu}
\author{Thang TQ L\^e}
\address{School of Mathematics \\
         Georgia Institute of Technology \\
         Atlanta, GA 30332-0160, USA}
\email{letu@math.gatech.edu}

\thanks{The authors were supported in part by National Science Foundation. \\
\newline
1991 {\em Mathematics Classification.} Primary 57N10. Secondary 57M25.
\newline
{\em Key words and phrases: hyperbolic volume conjecture, colored Jones
function, Jones polynomial, cyclotomic expansion, loop expansion,
WKB, $q$-difference equations, asymptotics, perturbation theory, 
Kontsevich integral.
}
}

\date{March 25, 2005 \hspace{0.5cm} First edition: March 25, 2005.}


\begin{abstract}
To a knot in 3-space, one can associate a sequence of Laurent polynomials,
whose $n$th term is the $n$th colored Jones polynomial. The 
Volume Conjecture for small angles states that the value of the 
$n$-th colored Jones polynomial at $e^{\a/n}$ is a sequence of 
complex numbers that grows subexponentially, for a fixed small complex angle 
$\a$. In an earlier publication, the authors
proved the Volume Conjecture for small purely imaginary angles,
using estimates of the cyclotomic expansion of a knot.
The goal of the present paper is to identify the polynomial growth
rate of the above sequence to all orders 
with the loop expansion of the colored Jones function.
Among other things, this provides a strong analytic form of the 
Melvin-Morton-Rozansky conjecture.
\end{abstract}

\maketitle

\tableofcontents


\section{Introduction}
\lbl{sec.intro}

\subsection{The volume conjecture for small angles}
\lbl{sub.volume}

In an earlier publication, the authors stated and proved the {\em Volume
Conjecture} 
for small purely imaginary angles; see \cite{GL2}. More precisely, the authors 
proved that for every knot $K$ in $S^3$ there exists a positive angle
$\a(K) >0$ such that 
\begin{equation}
\lbl{eq.VC}
\lim_{n \to \infty} \frac{\log|J_{K,n}(e^{\a/n})|}{n}=0
\end{equation}
for all $\a \in i [0, \a(K))$, where 
\begin{itemize}
\item
$f(e^{\a/n})$ denotes the evaluation of a rational function $f(q)$ at
$q=e^{\a/n}$, 
\item
$J_{K,n}(q) \in \BZ[q^{\pm}]$ is the {\em Jones polynomial} of a knot 
{\em colored} with the 
$n$-dimensional irreducible representation of $\mathfrak{sl}_2$,
normalized so that it equals to $1$ for the unknot (see \cite{J, Tu}).
\end{itemize}

In the following, we will refer to the complex parameter $\a$ as 
{\em the angle}, making contact with standard terminology from
hyperbolic geometry.
As was explained in \cite{GL2}, the above result agrees with the fact that
$$
\vol(\rho_{\a})=0
$$
where 
\begin{equation}
\lbl{eq.rho}
\rho_{\a}: \pi_1(S^3-K) \longto \SL_2(\BC),
\qquad \rho_{\a}(\mathfrak{m})= \mat {e^{\a/n}} 0 0
{e^{-\a/n}}. 
\end{equation}
is a {\em reducible} representation of the knot group in $\SL_2(\BC)$
with prescribed behavior on a meridian $\mathfrak{m}$ of the knot $K$.

For further reading concerning the history of the volume conjecture, we refer
the reader to \cite{Gu,K,MM}, as well as \cite{GL2}.

Notice that $\rho_{\a}$ is a 1-parameter deformation of the {\em trivial
representation} $\rho_0=I$.  

Moreover, Equation \eqref{eq.VC} implies that the sequence
$J_{K,n}(e^{\a/n})$ grows at a subexponential rate,
as $n$ approaches infinity, and $\a$ is small and purely imaginary.

The purpose of the present paper is to identify the polynomial
growth rate of $J_{K,n}(e^{\a/n})$ in terms of the 
inverse Alexander polynomial $\D_K$ of $K$, symmetrized by $\D_K(t^{-1})=
\D_K(t)$, and normalized by $\D_K(1)=1$, and $\D_{\text{unknot}}(t)=1$.
More precisely, we have the following theorem.

\begin{theorem}
\lbl{thm.11}
For every knot $K$ there exists an open neighborhood $U_K$ of $0 \in \BC$
such that for all complex angles $\a \in U_K$, we have:
\begin{equation}
\lbl{eq.thm11}
\lim_{n\to\infty}
J_{K,n}(e^{\a/n}) =
\frac{1}{\D_K(e^{\a})} \in \BC.
\end{equation}
Moreover, the convergence with respect to $\a$ is uniform on compact subsets
of $U_K$.
\newline
In particular since $\D_K(1)=1$, \eqref{eq.thm11} implies \eqref{eq.VC}.
\end{theorem}

The reader may compare the above theorem with the famous 
Melvin-Morton-Rozansky (MMR, in short) Conjecture, which was settled
by Bar-Natan and the first author in \cite{B-NG}. Let $\BQ[[h]]$ denote
the ring of {\em formal power series} in a variable $h$ with rational
coefficients.

\begin{theorem}
\lbl{thm.MMR}\cite{B-NG}
For every knot $K$ we have the following equality in the 
ring $\BQ[[h]]$:
\begin{equation}
\lbl{eq.MMR}
\lim_{n \to \infty} \,\, J_{K,n}(e^{h/n}) = \frac{1}{\D_K(e^{h})} 
\in \BQ[[h]],
\end{equation}
\end{theorem}

To avoid confusion, let us point out that Equation 
\eqref{eq.MMR} is a statement about coefficients of formal power series. 
In other words, \eqref{eq.MMR} can be phrased as follows:
for every $m \geq 0$, we have:
\begin{equation}
\lbl{eq.MMRalt}
\lim_{n \to \infty}
\coeff \left( J_{K,n}(e^{h/n}) , h^m \right)=
\coeff \left( \frac{1}{\D_K(e^{h})} , h^m\right),
\end{equation}
where for an analytic function $f(x)$ we define:
$$
\coeff( f(h), h^m)=\frac{1}{m!} \frac{d^m}{dh^m}|_{h=0} f(h).
$$
Actually, for every $m \geq 0$, $
\coeff \left( J_{K,n}(e^{h/n}) , h^m \right) $ is a polynomial in $1/n$
of degree $m$ (see also Section \ref{sub.fti} below). Thus, the limit
with respect to $n \to\infty$ in \eqref{eq.MMRalt} exists and is simply
the constant term of the above-mentioned polynomial. Identifying that constant
term with the right hand side of \eqref{eq.MMRalt} is the non-trivial
part of the MMR Conjecture.
 
Let us compare Theorems \ref{thm.11} and \ref{thm.MMR}.
Since convergence with respect to $\a$ is uniform on compact subsets,
it is easy to see that Theorem \ref{thm.11} implies Theorem \ref{thm.MMR}.
In that sense, we may say
that Theorem \ref{thm.11} is an analytic form of the MMR Conjecture.

Thus, Theorem \ref{thm.11} can be viewed as a statement about the volume 
conjecture for small angles, as well as an analytic form of the MMR Conjecture.

Armed with Theorem \ref{thm.11} one may ask for a full asymptotic expansion
of the left hand side of \eqref{eq.thm11} in terms of powers of $1/n$.
Before we answer this question, let us recall what is known on the level
of formal power series, that is, about the $1/n$ terms of 
\eqref{eq.MMRalt}.

Rozansky discovered that after resummation, for every fixed $m \geq 0$,
the $1/n^m$ terms of \eqref{eq.MMRalt} are rational
functions in a variable $e^h$. Let us state Rozansky's discovery concretely.

\begin{theorem}
\lbl{thm.Zrat}\cite{Ro}
For every knot $K$ there exists a sequence  $P_{K,k}(q) \in \BQ[q^{\pm}]$
of Laurent polynomials with $P_{K,0}(q)=1$ such that 
\begin{equation}
\lbl{eq.Zrat}
 J_{K,n}(e^{h/n}) \sim_{n \to \infty} 
\sum_{k=0}^\infty 
\frac{P_{K,k}(e^{h})}{\D_K(e^{h})^{2k+1}} \left(\frac{h}{n}\right)^k 
\in \BQ[[h]]
\end{equation}
in the ring $\BQ[[h]]$ of formal power series in $h$.
\end{theorem}

A different proof, valid for all simple Lie groups, was given in \cite{Ga1},
using work of \cite{GK}.

Let us point out that \eqref{eq.Zrat} means the following: 
for every $N \geq 0$ we have:
\begin{equation}
\lbl{eq.Zratalt}
\lim_{n\to\infty}
\left( \frac{n}{h} \right)^N \left(
J_{K,n}(e^{h/n})-\sum_{k=0}^{N-1}
\frac{P_{K,k}(e^{h})}{\D_K(e^{h})^{2k+1}} \left(\frac{h}{n}\right)^k
\right)=
\frac{P_{K,N}(e^h)}{\D_K^{2N+1}(e^h)} \in \BQ[[h]].
\end{equation}

\subsection{Asymptotics to all orders}
\lbl{sub.results}

Our results are the following:

\begin{theorem}
\lbl{thm.1}
For every knot $K$ there exists an open neighborhood $U_K$ of $0 \in \BC$
such that for all complex angles $\a \in U_K$, we have
an asymptotic expansion (uniform on compact subsets of $U_K$
with respect to $\a$):
\begin{equation}
\lbl{eq.thm1}
J_{K,n}(e^{\a/n}) \sim_{n \to \infty} \sum_{k=0}^\infty 
\frac{P_{K,k}(e^{\a})}{\D_K(e^{\a})^{2k+1}} \left(\frac{\a}{n}\right)^k
.
\end{equation}
\end{theorem}

In other words, for $\a \in U_K$ and every $N \geq 0$,
\begin{equation}
\lbl{eq.thm1alt}
\lim_{n\to\infty}
\left( \frac{n}{\a} \right)^N \left(
J_{K,n}(e^{\a/n})-\sum_{k=0}^{N-1}
\frac{P_{K,k}(e^{\a})}{\D_K(e^{\a})^{2k+1}} \left(\frac{\a}{n}\right)^k
\right)=
\frac{P_{K,N}(e^{\a})}{\D_K^{2N+1}(e^{\a})} \in \BC.
\end{equation}
Moreover,  convergence with respect to $\a$ is uniform on compact subsets
of $U_K$.

Thus, the above theorem determines to all orders the asymptotic expansion
of the volume conjecture for small angles.

\subsection{A small dose of physics}
\lbl{sub.physics}

One does not need to know the relation of the colored Jones
function and quantum field theory in order to understand the statement
and proof of Theorem \ref{thm.1}. Nevertheless, we want to add some 
philosophical comments, for the benefit of the willing reader. 
According to Witten (see \cite{Wi}), the Jones polynomial
$J_{K,n}$ can be expressed by a partition function of a topological quantum
field theory in $3$ dimensions---a gauge theory with Chern-Simons Lagrangian.
The stationary points of the Lagrangian
correspond to $SU(2)$-flat connections on an ambient manifold, and the
observables are knots, colored by the $n$-dimensional irreducible 
representation of $SU(2)$. In case of a knot in $S^3$, there is only one
ambient flat connection, and the corresponding perturbation theory
is a formal power series in $h=\log q$.

Rozansky exploited a cut-and-paste property of the Chern-Simons path integral
and considered perturbation theory of the knot complement, along an abelian
flat connection with monodromy given by \eqref{eq.rho}. In fact, Rozansky
calls such an expansion the $U(1)$-{\em RCC connection contribution} to the
Chern-Simons path integral, where RCC stands for {\em reducible connection
contribution}, and $U(1)$ stands for the fact that the flat $SU(2)$
connections are actually $U(1)$-valued abelian connections.
Formal properties of such a perturbative expansion, enabled Rozansky to
deduce (in physics terms) the loop expansion of the colored Jones function.
In a later publication, Rozansky proved the existence of the loop expansion
using an explicit state-sum description of the colored Jones function.

Of course, perturbation theory means studying formal power series that rarely
converge. Perturbation theory at the trivial flat connection in a knot
complement converges, as it resums to a Laurent polynomial in $e^h$; namely
the colored Jones polynomial. 
The volume conjecture for small complex angles is precisely the
statement that perturbation theory for abelian flat connections (near the 
trivial one) does converge.

At the moment, there is no physics (or otherwise) formulation of perturbation
theory of the Chern-Simons path integral along a discrete and faithful
$\SL_2(\BC)$ representation. Nor is there an adequate explanation of the
relation between $SU(2)$ gauge theory (valid near $\a=0$) and 
a complexified $\SL_2(\BC)$ gauge theory, valid near $\a=2\pi i$.
These are important and tantalizing questions, with no answers at present.

\subsection{WKB}
\lbl{sec.WKB}

Since we are discussing physics interpretations of Theorem \ref{thm.1}
let us make some more comments. Obviously, when the angle
$\a$ is sufficiently big, the asymptotic expansion of Equation \eqref{eq.thm1}
may break down. For example, when $e^{\a}$ is a complex root of the Alexander 
polynomial, then the right hand side of \eqref{eq.thm1} does not make sense,
even to leading order. In fact, when $\a$ is near $2 \pi i$, then the solutions
are expected to grow exponentially, and not polynomially, according to the
Volume Conjecture.

The breakdown and change of rate of asymptotics is 
a well-documented phenomenon well-known in physics, associated with WKB
analysis, after Wentzel-Krammer-Brillouin; see for example \cite{O}.
 In fact, one may obtain an independent proof of Theorem \ref{thm.1}
using {\em WKB analysis}, that is, the study of asymptotics of solutions of 
difference equations with a small parameter. 
The key idea is that the sequence of colored Jones  
functions is a solution of a linear $q$-difference equation, as was 
established in \cite{GL1}. A discussion on WKB analysis of $q$-difference
equations was given by Geronimo and the first author in \cite{GG}.

The WKB analysis can, in particular, determine {\em small exponential
corrections} of the form $e^{-c_{\a} n}$ to the asymptotic expansion of 
Theorem \ref{thm.1}, where $c_{\a}$ depends on $\a$, with $\text{Re}(c_{\a})
<0$ for $\a$ sufficiently small. 
These exciting small exponential corrections
cannot be captured by classical asymptotic analysis (since they vanish
to all orders in $n$), but they are important and dominant (i.e.,
$\text{Re}(c_{\a})>0)$ when $\a$ is
near $2 \pi i$, according to the volume conjecture. Understanding the change of
sign of $\text{Re}(c_{\a})$ past certain so-called Stokes directions 
is an important question that WKB addresses.

We will not elaborate or use the WKB analysis in the present paper.
Let us only mention that the loop expansion of the colored Jones function
can be interpreted as WKB asymptotics on a $q$-difference equation satisfied
by the colored Jones function.

\subsection{The main ideas}
\lbl{sub.main}

The main ideas of Theorem \ref{thm.1} is to compare three different
views of the Jones polynomial: one coming from perturbative quantum field
theory, one from a resummation of quantum field theory (known as the loop
expansion), and a third non-perturbative view, in terms of the cyclotomic
function.

The main advantage of the cyclotomic function of a knot is a key integrality
property, due to Habiro, and a priori exponential estimates for the $l^1$-norm
and quadratic bounds for the degrees of the revelant polynomials. 
The latter were established in \cite{GL2}. Using these bounds, we can prove 
that for small enough complex angles, a sequence of holomorphic functions
is uniformly bounded, and the limit of derivatives of any order (at zero)
exists; see Theorem \ref{thm.boundC3}.
A key lemma from complex analysis on normal families  
guarantees under the above hypothesis that the sequence of holomorphic 
functions converges, uniformly on compact sets, to a holomorphic function
whose derivatives (at zero) are the limits of the derivatives of the original
sequence of holomorphic functions.

\subsection{Acknowledgement}
Soon after the completion of the authors' work \cite{GL2}, H. Murakami
posted an interesting paper, in which he identified the polynomial growth
of the volume conjecture for small angles, for the case of the $4_1$ knot; 
see \cite{M}.
Upon reading Murakami's paper, it became clear that the methods of \cite{GL2}
can be adapted to all knots, and to all orders, for small complex angles.
We wish to thank Murakami who motivated our present work.

\section{Three expansions of the Jones polynomial}
\lbl{sec.fti}

\subsection{Finite type invariants and the Jones polynomial}
\lbl{sub.fti}

The colored Jones function of a knot is a 2-parameter invariant, that
depends on the color $n$ and the formal parameter 
$$
h=\log q.
$$ 
{\em Perturbative quantum field theory} (formalized mathematically by the
{\em Kontsevich integral} of a knot, and its image under the $\mathfrak{sl}_2$
{\em weight system}, described for example in \cite{B-N}) 
gives the following expansion of the colored Jones function:

\begin{eqnarray}
\lbl{eq.pert}
J_{K,n}(e^h) &=& \sum_{0 \leq i, 0 \leq j \leq i} c_{K,i,j} n^j h^i \\
\notag
&=& \sum_{0 \leq i, 0 \leq j \leq i} c_{K,i,j} (nh)^j h^{i-j} \\
\notag
&=& \sum_{0 \leq j, k} c_{K,j+k,j} (nh)^j h^k.
\end{eqnarray}
Here, $K\to c_{K,i,j}$ are {\em finite type knot invariants} of type $i$;
see \cite{B-N}.  
The important property is that $c_{K,i,j}=0$ in the $(i,j)$ plane and
above the diagonal $i=j$. Thus, one can resum the formal power series 
as follows:
\begin{eqnarray}
\lbl{eq.R}
J_{K,n}(e^h) &=& \sum_{k=0}^\infty R_{K,k}(nh) h^k,
\end{eqnarray}
where
$$
R_{K,k}(x)=\sum_{0 \leq j} c_{K,j+k,j} x^j \in \BQ[[x]].
$$

\subsection{The loop expansion of the Jones polynomial}
\lbl{sub.loop}

The Melvin-Morton-Rozansky Conjecture states that
$$
R_{K,0}(x)=\frac{1}{\D_K(e^x)}.
$$
More generally, in \cite{Ro}, Rozansky proves that
$$
R_{K,k}(x)=\frac{P_{K,k}(e^x)}{\D_K(e^x)^{2k+1}}
$$
for Laurent polynomials $P_{K,k}(q) \in \BQ[q^{\pm}]$.

Although the polynomials $P_{K,k}(q)$ are not finite type invariants
(with respect to the usual crossing change of knots), they are indeed
finite type invariants with respect to a loop move described in \cite{GR}.
We will not use this fact in our paper.

Rozansky conjectured that the resummation given by the above equations
could be preformed on the level of a universal perturbative invariant
(the Kontsevich integral of a knot; see \cite{B-N}), and this was proven
to be the case in \cite{GK}. As a result, one obtains a proof of this 
resummation property valid for all simple Lie algebras, see \cite{Ga1}.

\subsection{The cyclotomic expansion of the Jones polynomial}
\lbl{sub.cyclotomic}

In \cite{Ha}, Habiro introduced an alternative packaging of the colored Jones
function $J_{K,n}$; using the so-called {\em cyclotomic function} $C_{K,n}$.
The latter is related to the former by the following 

\begin{equation}
\lbl{eq.J2C}
J_{K,n}(q)=\sum_{k=0}^n C_{n,k}(q) C_{K,k}(q),
\end{equation}
where
\begin{eqnarray*}
\lbl{eq.cyclokernel}
C_{n,k}(q) &:=& \frac{1}{q^{n/2}-q^{-n/2}}
\prod_{j=n-k}^{n+k} (q^{j/2}-q^{-j/2}) \\ 
& = & \prod_{j=1}^k (( q^{n/2}-q^{-n/2})^2 - (q^{j/2}-q^{-j/2})^2) \\
& = & \prod_{j=1}^k (( q^{n/2}+q^{-n/2})^2 - (q^{j/2}+q^{-j/2})^2). 
\end{eqnarray*}

Thus, in a sense $J_{K,n}$ and $C_{K,n}$ are related by a lower-diagonal
invertible matrix. For an explicit inversion of the above equation (which we
will not use in the present paper), we refer the reader to \cite[Sec.4]{GL1}.

\subsection{Comparing the cyclotomic and the loop expansion}
\lbl{sub.comparing}

So far, we have three expansions: the finite type expansion, the loop 
expansion and the cyclotomic expansion. Now, we'll compare the last two.
In other words, we'll compare Equations \eqref{eq.R} and \eqref{eq.J2C}.

Let 
$$
q=e^h, \qquad x=nh.
$$ 
For a function $f(q)$, let us denote by $\la f \ra_k$ the $k$-th coefficient
in the Taylor expansion of $f(e^h)$ around $h=0$. 
Of course,
$$
\la f \ra_k=\frac{1}{k!} \frac{d^k}{dh^k}|_{h=0} f(e^h).
$$
In other words, we have:
$$
f(e^h)=\sum_{k=0}^\infty \la f \ra_k h^k \in \BQ[[h]].
$$

\begin{lemma}
\lbl{lem.compare2}
\rm{(a)} For every knot $K$, we have the following equality
in $\BQ[[x,h]]$:
$$
\sum_{k=0}^\infty R_{K,k}(x) h^k = \sum_{k=0}^\infty C_{K,k}(e^h)
\prod_{j=1}^k (e^{x/2}-e^{-x/2})^2-(e^{jh/2}-e^{-jh/2})^2) \in \BQ[[x,h]].
$$
\rm{(b)} It follows that for every $k$,
$$
R_{K,k}(x)=\sum_{l=0}^\infty 
\sum_{j=0}^k \la C_{K,l} \ra_j  z^{2l-[j/2]} p_{l,j,k}(z)
$$
where
$$
z=e^{x/2}-e^{-x/2},
$$
and $p_{l,j,k}(z)$ is an even polynomial of $z$ of degree $[j/2]$, with
coefficients polynomials of $l$ of degree $k+1$.
\newline
\rm{(c)} 
In particular, we have:
\begin{eqnarray*}
R_{K,0}(x) &=& \sum_{l=0}^\infty \la C_{K,l} \ra_0 z^{2l} \\
R_{K,1}(x) &=& \sum_{l=0}^\infty \la C_{K,l} \ra_1 z^{2l} \\
R_{K,2}(x) &=& \sum_{l=0}^\infty \la C_{K,l} \ra_2 z^{2l} -
\sum_{l=0}^\infty \la C_{K,l} \ra_0 \frac{l(l+1)(2l+1)}{6} z^{2l-2} \\
R_{K,3}(x) &=& \sum_{l=0}^\infty \la C_{K,l} \ra_3 z^{2l} -
\sum_{l=0}^\infty \la C_{K,l} \ra_1 \frac{l(l+1)(2l+1)}{6} z^{2l-2}
\end{eqnarray*}
in $\BQ[[x]]$.
\end{lemma}

\begin{proof}
It follows easily, working in the ring $\BQ[[x,h]]$, and using the fact
that the map:
$$
\BQ(e^x)[[h]] \longto \BQ[[x,h]]
$$
given by $e^x=\sum_{k=0}^\infty x^k/k!$ is 1-1.
\end{proof}

\section{Proof of Theorem \ref{thm.11}}
\lbl{sec.proofs}

Let us assume for the moment the following theorem, whose proof
will be given in the next section.

\begin{theorem}
\lbl{thm.boundC3}
\rm{(a)}
For every knot $K$ there exist an open neighborhood $U_K$ of $0 \in \BC$ 
and a positive number $M$ such
that for $\a \in U_K$, and all $n \geq 0$, we have:
$$
|J_{K,n}(e^{\a/n})| < M.
$$
\rm{(b)}
Moreover, for every $m \geq 0$, the following limit exists and given by:
$$
\lim_{n\to\infty} \frac{d^m}{d \a^m}|_{\a=0} J_{K,n}(e^{\a/n})
=m! \,\, \coeff\left( \frac{1}{\D_K(e^{\a})}, \a^m \right).
$$
\end{theorem}

\subsection{A lemma from complex analysis}
\lbl{sub.complex}

The proof of Theorem \ref{thm.1} will use the following lemma on normal 
families that is sometimes refered to by the name of Vitali and 
Montel's theorem. For a reference, see \cite{Hi,Sch}. The lemma exhibits
the power of holomorphy, coupled with uniform boundedness.

Let $\D_r=\{z \in \BC \, : \, |z| < r \}$ 
denote the open complex disk around $0$ of radius $r >0$.

\begin{lemma}
\lbl{lem.complex}
If $f_n: \D_r \to \bar\D_M$ is a sequence of holomorphic functions such
that for every $m \geq 0$, we have:
$$
\lim_{n \to \infty} f^{(m)}_n(0) =a_m.
$$
Then, 
\begin{itemize}
\item
The limit $f(z)=\lim_n f_n(z)$ exists pointwise for $z \in D_r$.
\item
$f: D_r \to \bar\D_M$ is holomorphic,
\item
The convergence is uniform on compact subsets, and
\item
For every $m$, $f^{(m)}(0)=a_m$.
\end{itemize}
\end{lemma}

\begin{proof}
$\{f_n\}_n$ is uniformly bounded, so it is a normal family, and contains
a convergent subsequence $f_j\to f$. Convergence is uniform on compact sets,
and $f$ is holomorphic, and for every $m \geq 0$, 
$\lim_j f_j^{(m)}(0)=f^{(m)}(0)=a_m$. 

If $\{f_n\}_n$ is not convergent, since it is a normal family,
then there exist two subsequences that converge to $f$ and $g$ respectively,
with $f \neq g$. 
Applying the above discussion, it follows that $f$ and $g$ are holomorphic
functions with equal derivatives of all orders at $0$. Thus, $f=g$, giving
a contradiction. Thus, $\{f_n\}_n$ is convergent and the result follows 
from the above discussion.
\end{proof}

\begin{remark}
\lbl{rem.necessarynormal}
We have seen that the hypotheses in Lemma \ref{lem.complex} are sufficient
to ensure existence of the limit and uniform convergence on compact sets.
It is easy to see that these hypotheses are also necessary.
\end{remark}

\subsection{Proof of Theorem \ref{thm.11}}
\lbl{sub.proofthm11}

Fix a knot $K$ and an open neighborhood $U_K$ of $0 \in \BC$ as in 
Theorem \ref{thm.boundC3}. 
Theorem \ref{thm.boundC3} and Lemma \ref{lem.complex} imply that
for $\a \in U_K$, 
$$
\lim_{n\to\infty}J_{K,n}(e^{\a/n}) 
= \frac{1}{\D_K(e^{\a})}.
$$
Moreover, convergence with respect to $\a$ is uniform on compact subsets
of $U_K$. This proves Theorem \ref{thm.11}.
\qed

\section{Estimates of the cyclotomic function}
\lbl{sec.estimates}

This section is devoted to the proof of Theorem \ref{thm.boundC3}.
Our main tool will be estimates in the cyclotomic expansion of a knot,
similar to the ones used in \cite{GL2}.

A key result of Habiro is an {\em integrality property} 
of the cyclotomic function $n\to C_{K,n}$ of a knot. Namely, 
$$
C_{K,n}(q) \in \BZ[q^{\pm}]
$$ 
for all knots $K$ and all $n$; see \cite{Ha}. 

We will use two further results from \cite{GL2}: an exponential 
bound on the size of the coefficients of $C_{K,n}$, 
and a quadratic bound on the min and max degrees of $C_{K,n}$. 
Recall that for a Laurent polynomial 
$f(q)=\sum_k a_k q^k$, we define its $l^1$ norm by 
$$
||f||_1=\sum_k |a_k|.
$$

\begin{theorem}
\lbl{thm.boundC}
\rm{(a)} For every knot $K$ we have:
\begin{equation}
\lbl{eq.LC}
||C_{K,n}||_1 \leq e^{C n + C' \log n}
\end{equation}
\rm{(b)} Moreover, 
$$
\mathrm{maxdeg}_q (C_{K,n}) =O(n^2), \qquad \mathrm{mindeg}_q (C_{K,n}) 
=O(n^2).
$$
\end{theorem}

Here, and below, the $O(f(n))$ notation means that a quantity bounded
by a constant times $f(n)$.

\begin{theorem}
\lbl{thm.boundC2}
For every knot $K$,  there exist constants $C, C', C''$
and $C'''$ (that depend on $K$) such that for all $n \geq 0$ and $k \geq 0$ 
we have:
\begin{equation}
\lbl{eq.boundC2}
|C_{K,n}^{(k)}(e^{\a})| \leq e^{C n + C' (k+1) 
\log n + \Re(\a) C'' n^2 + C'''},
\end{equation}  
where $C_{K,n}^{(k)}$ denotes the $k$-th derivative of $C_{K,n}(e^h)$ 
with respect to $h$.
\end{theorem}

\begin{proof}
Let us write 
$$
C_{K,n}(q)=\sum_{j=-C_1 n^2}^{C_1 n^2} a_{j,n} q^j.
$$
Then, 
$$
C_K^{(k)}(e^h)=\sum_{j=-C_1 n^2}^{C_1 n^2} a_{j,n} j^k e^{jh}.
$$
We will estimate each coefficient and each monomial by:
\begin{eqnarray*}
|a_{j,n}| & \leq & ||C_{K,n}||_1 \leq e^{C n + C' \log n} \\
|j|^k & \leq & (C_1 n^2)^k \\
|e^{\a j}| & \leq & e^{|\Re(\a)| C_1 n^2}.
\end{eqnarray*}
The result follows.
\end{proof}

\begin{corollary}
\lbl{cor.boundC2}
With the notation of Theorem \ref{thm.boundC2}, for every $n \geq 0$
and $0 \leq k \leq n$, and $0 \leq l \leq n$, we have:
$$
|C_{K,k}^{(l)}(e^{\a/n})| \leq e^{C k + C'(l+1) \log k 
+ |\Re(\a)| C'' k + C'''}
$$
\end{corollary}

Let us recall an elementary estimate from \cite[Sec.3]{GL2}.

\begin{lemma}
\lbl{lem.estimate}
There exist positive constants $C_1, C_2$ and $C_3$, so that for all complex
numbers $\a$ with $0 < \Re(\a) < 1/6$, and for every $0 \leq k < n$
we have: 
$$
|C_{n,k}(e^{\a/n})| \leq e^{C_1 k \log|\a| + C_2 \log k + C_3}.
$$
\end{lemma}

\begin{proof}(of Theorem \ref{thm.boundC3})
Combining Corollary \ref{cor.boundC2}
and Lemma \ref{lem.estimate}, it follows that for all $0 \leq k \leq n$,
we have:
$$
|C_{n,k}(e^{\a/n}) C_{K,k}(e^{\a/n})| \leq 
e^{C k + C' \log k + |\Re(\a)| C'' k + C''' + C_1 k \log|\a| +
C_2 \log k + C_3}.
$$
Let us choose $\a \in U_K$, where
\begin{equation}
\lbl{eq.UK}
U_K = \{\a \in \BC \, | \, C+C'' |\Re(\a)| + C_1 \log|\a| <0 \}.
\end{equation}
Then, equation \eqref{eq.J2C} and the
above estimate conclude the first part of Theorem \ref{thm.boundC3}.

The second part follows from Equation \eqref{eq.R} and the MMR Conjecture. 
Indeed, consider the sequence 
$$
f_n: U_K \to \{z: \, |z| < N\}, \qquad 
\a \to f_n(\a)=J_{K,n}(e^{\a/n}).
$$
Since $J_{K,n}(q)$ is a Laurent polynomial in $q$, it follows that 
$f_n$ is an entire function. Equation \eqref{eq.R} implies that
$$
f_n(\a)=\sum_{k=0}^\infty R_{K,k}(\a) \left( \frac{\a}{n} \right)^k.
$$
Thus, for every $m \geq 0$, 
$$
f_n^{(m)}(0)=m! \left( \coeff(R_{K,0}(\a),\a^m) + 
\frac{1}{n} \coeff(R_{K,1}(\a),\a^{m-1}) + \dots 
\frac{1}{n^m} \coeff(R_{K,m}(\a),\a^0) \right).
$$
Thus, using the MMR Conjecture, we obtain:
\begin{eqnarray*}
\lim_{m\to\infty} f_n^{(m)}(0) &=& m! \,\, \coeff(R_{K,0}(\a),\a^m) \\
&=& m! \,\,  \coeff \left( \frac{1}{\D_K(e^{\a})},\a^m \right).
\end{eqnarray*}
The result follows.
\end{proof}

\section{Proof of Theorem \ref{thm.1}}
\lbl{sec.allorders}

To leading order (i.e., $N=0$ in \eqref{eq.Zratalt}) Theorem \ref{thm.1}
is Theorem \ref{thm.11}. 
By now, it should be clear the strategy for proving Theorem \ref{thm.1}
to all orders. To simplify notation, let us define:
\begin{equation}
\lbl{eq.JN}
J_{K,n}^{(N)}(e^{\a/n})=J_{K,n}(e^{\a/n})-\sum_{k=0}^{N-1}
\frac{P_{K,k}(e^{\a})}{\D_K(e^{\a})^{2k+1}} \left(\frac{\a}{n}\right)^k.
\end{equation}
Theorem \ref{thm.1} follows from the following result
and the argument of Section \ref{sub.proofthm11}.

\begin{theorem}
\lbl{thm.boundC4}
\rm{(a)}
For every knot $K$ 
there exists an open neighborhood $U_K$ of $0 \in \BC$ such that
for every $N \geq 0$ there exists a positive number $M_N$ such
that for $\a \in U_K$, and all $n \geq 0$, we have:
$$
\left|\left( \frac{n}{\a} \right)^N 
J_{K,n}^{(N)}(e^{\a/n}) \right| < M_N.
$$
\rm{(b)}
Moreover, for every $m \geq 0$, the following limit exists and given by:
$$
\lim_{n\to\infty} \frac{d^m}{d \a^m}|_{\a=0} 
\left( \left( \frac{n}{\a} \right)^N J_{K,n}^{(N)}(e^{\a/n}) \right)
=m! \,\, \coeff\left( 
\frac{P_{K,N}(e^{\a})}{\D_K(e^{\a})^{2N+1}}, \a^m \right).
$$
\end{theorem}

\begin{proof}
We will prove the theorem by induction on $N$. For $N=0$, this
is Theorem \ref{thm.11} proven in Section \ref{sec.proofs}. 
Let us assume that it is true for $N-1$. 

Let us define for every $k \geq 0$,
two auxiliary biholomorphic functions
\begin{eqnarray*}
c_k(x,\e) &=&
\prod_{j=1}^k (e^{x/2}-e^{-x/2})^2-(e^{jh/2}-e^{-jh/2})^2), \\
g_{K,k}(x,\e) &=& c_k(x,\e) C_{K,k}(e^{\e}).
\end{eqnarray*}
Thus, using the definition of $C_{n,k}$ and Equation \eqref{eq.J2C},
it follows that:
\begin{equation}
\lbl{eq.now1}
C_{n,k}(e^{\a/n})=c_k(\a,\a/n), \qquad 
J_{K,n}(e^{\a/n})=\sum_{k=0}^n g_{K,k}(\a,\a/n).
\end{equation}
For a function $h=h(x)$, let us define the $N$-th Taylor approximation by:
$$
\Taylor^N(h,x)=\sum_{j=0}^N \frac{h^{(j)}(0)}{j!} x^j.
$$
Applying Lemma \ref{lem.compare2} to the function $\e\to g_{K,k}(\a,\e)$,
and evaluating at $\e=\a/n$, it follows that:
\begin{eqnarray}
\sum_{k=0}^{N-1}
\frac{P_{K,k}(e^{\a})}{\D_K(e^{\a})^{2k+1}} \left(\frac{\a}{n}\right)^k &=&
\sum_{k=0}^\infty \Taylor^{N-1}(g_{K,k}(\a,\cdot), \frac{\a}{n}) \\
&=&
\lbl{eq.now2}
\sum_{k=0}^n \Taylor^{N-1}(g_{K,k}(\a,\cdot), \frac{\a}{n}) 
+ \text{err}_n(\a).
\end{eqnarray}
Equations \eqref{eq.JN}, \eqref{eq.now1} and \eqref{eq.now2} and Taylor's 
theorem imply that:
\begin{eqnarray*}
J_{K,n}^{(N)}(e^{\a/n})&=&
J_{K,n}(e^{\a/n})-\sum_{k=0}^{N-1}
\frac{P_{K,k}(e^{\a})}{\D_K(e^{\a})^{2k+1}} \left(\frac{\a}{n}\right)^k \\
&=& \sum_{k=0}^n g_{K,k}(\a,\a/n)- 
\sum_{k=0}^n \Taylor^{N-1}(g_{K,k}(\a,\cdot), \frac{\a}{n}) 
- \text{err}_n(\a) \\
&\approx& \left(\frac{\a}{n}\right)^N \sum_{k=0}^n 
\,
\frac{1}{N!} \frac{\pt^N}{\pt \e^N}|_{\e \approx \a/n}
g_{K,k}(\a,\e)
-\text{err}_n.
\end{eqnarray*}

The analytiticy of $g_{K,k}$ and Theorem \ref{thm.boundC2} implies that
there exists a positive $M'_N$ such that for all $\a \in U_K$
(defined in \ref{eq.UK}), we have:
$$
|\text{err}_n(\a)| < M'_N.
$$
Corollary \ref{cor.boundC2} and Equation \eqref{eq.now1} imply that there 
exists a positive $M_N$ such that 
$$
|\left(\frac{n}{\a}\right)^N J_{K,n}^{(N)}(e^{\a/n})| < M_N
$$
for all $n \geq 0$ and for all $\a \in U_K$. 
This proves part (a) of Theorem \ref{thm.boundC4}.

For part (b), we will use Equation \eqref{eq.R}, which implies that:
$$
J_{K,n}^{(N)}(e^{\a/n})=\sum_{k=N}^\infty R_{K,k}(\a) 
\left(\frac{\a}{n}\right)^{k}.
$$
Thus, for every $m \geq 0$, 
\begin{eqnarray*}
\frac{d^m}{d \a^m}|_{\a=0} 
\left( \left( \frac{n}{\a} \right)^N J_{K,n}^{(N)}(e^{\a/n}) \right)
&=& m! \left( \coeff(R_{K,N}(\a),\a^m) + 
\frac{1}{n} \coeff(R_{K,N+1}(\a),\a^{m-1}) + \right. \\ & & \dots +
\left. \frac{1}{n^m} \coeff(R_{K,N+m}(\a),\a^0) \right).
\end{eqnarray*}
Using Rozansky's theorem \ref{thm.Zrat} and Equation \eqref{eq.Zratalt},
we obtain:
\begin{eqnarray*}
\lim_{m\to\infty} \frac{d^m}{d \a^m}|_{\a=0} 
\left( \left( \frac{n}{\a} \right)^N J_{K,n}^{(N)}(e^{\a/n}) \right)
&=& m! \,\, \coeff (R_{K,N}(\a),\a^m) \\
&=&
m! \,\, \coeff \left(\frac{P_{K,N}(e^{\a})}{\D_K(e^{\a})^{2N+1}},
\a^m \right).
\end{eqnarray*}
The result follows.
\end{proof}

\ifx\undefined\bysame
    \newcommand{\bysame}{\leavevmode\hbox
to3em{\hrulefill}\,}
\fi

\end{document}